\documentclass[a4paper, 12pt]{amsart}
\usepackage{amsmath, amssymb}
\usepackage{color}
\usepackage{verbatim}
\usepackage{graphicx}
\usepackage{marvosym}

\newcommand{\R}{\mathbb{R}}

\newcommand{\swrz}{\mathcal{S}}

\newcommand{\field}[1]{\mathbb{#1}} 

\newcommand{\bZ}{\field{Z}} 
\newcommand{\bR}{\field{R}} 
\newcommand{\E}{\mathrm{E}} 
\renewcommand{\Re}{\mathrm{Re}} 
\renewcommand{\Im}{\mathrm{Im}} 

\def\figWidth{0.45\textwidth} 

\title{The pole behaviour of the phase derivative of the short-time Fourier transform}
\author{Peter Balazs$^{A)}$, Dominik Bayer$^{B)}$, Florent Jaillet$^{C)}$ and Peter S{\o}ndergaard$^{D)}$}
\thanks{$^{A)}$ Acoustics Research Institute, Austrian Academy of Sciences, Vienna A-1040, Austria; \\
$^{B)}$ Corresponding author; Acoustics Research Institute, Austrian Academy of Sciences, Vienna A-1040, Austria;\hspace{2ex}  \Telefon\ +43 1 51581-2517 \hspace{2ex} \Letter\ bayerd@kfs.oeaw.ac.at; \\
$^{C)}$ Institut de Neurosciences de la Timone UMR 7289, Aix Marseille Universit\'e, CNRS, 13385 cedex 5, Marseille,
France; \\ 
$^{D)}$ Oticon A/S, 2765 Sm{\o}rum, Denmark;
}

\newtheorem{theorem}{Theorem}[section]
\newtheorem{prop}[theorem]{Proposition}
\newtheorem{definition}[theorem]{Definition}
\newtheorem{corollary}[theorem]{Corollary}
\newtheorem{lemma}[theorem]{Lemma}

\begin{document}

\begin{abstract}
The short-time Fourier transform (STFT) is a time-frequency representation
widely used in applications, for example in audio signal processing. Recently it has been shown that not only the amplitude, but also the phase of this representation can be successfully exploited for improved analysis and processing. 
In this paper we describe a rather peculiar pole phenomenon in the phase derivative, a recurring pattern that appears in a characteristic way in the neighborhood around any of the zeros of the STFT, a negative peak followed by a positive one. We describe this phenomenon numerically and provide a complete analytical explanation.   
\end{abstract}

\maketitle

\section{Introduction}

\label{sec:Introduction}

The {\em short-time Fourier transform} (STFT) \cite{carmtorr1,fland99} is a
time-frequency representation widely used in audio signal
processing. A common definition of the STFT\footnote{This is the {\em
  frequency-invariant} STFT.} is \begin{flalign}
\label{eq:defV}
V(f,g)(x,\omega)= \int f(t) \overline{g(t-x)} e^{-2\pi i\omega t} \,dt.
\end{flalign}
The STFT $V(f,g)(x,\omega)$ provides information about the frequency
content of the signal $f$ at time $x$ and frequency $\omega$. The
analyzing window $g$ determines the resolution in time and frequency. 

The interpretation of the modulus of the STFT is relatively easy,
considering the fact that the spectrogram (defined as the square
absolute value of the STFT) can be interpreted as a time-frequency
distribution of the signal energy. This interpretation led to the
important success of the STFT in signal processing. In particular, it
has been widely used for applications in speech processing and
acoustics as a graphical tool for signal analysis \cite{quat01}.

But the interpretation of the phase of the STFT is less obvious, and
was thus hardly considered in applications for some time.

The phase can be
of particular interest for certain applications, as illustrated by
important applications such as phase vocoder \cite{phavoc66,dols1} or
reassignment \cite{kodera1976nmn,aufl95}.
In digital image processing it is well known that the phase
information of the discrete Fourier transform is at least as important
as the amplitude information.
In \cite{oppenheim81} it is shown that
as long as the phase of the discrete Fourier transform of an image is
retained and the amplitude is set to 1, the image can still be
recognized. Similar effects can also be shown for acoustic signal depending on the parameters of the STFT \cite{xxlwau2}.

For applications modifying the STFT coefficients, phase
information is essential again. For these types of applications, in
particular for the applications using Gabor frame multipliers
\cite{feinow1,xxlmult1} which motivated the present study, better
understanding of the structure of the phase is necessary to improve
the processing possibilities.

The phase of the STFT is usually not considered directly. In fact, it
is more interesting to consider the phase derivative over time or
frequency. Indeed, these quantities appear naturally in the context of
reassignment \cite{aufl95} and manipulations of phase derivative over
time is the idea behind the phase vocoder \cite{dols1}. Their
interpretation is easier, as the derivative of phase over time can be
interpreted as local instantaneous frequency while the derivative of
the phase over frequency can be interpreted as a local group delay.

To numerically compute the local instantaneous frequency, an unwrapping
of the phase is needed to avoid discontinuities. This is the classical
method used in \cite{dols1,kodera1976nmn}. Another method was found
 in \cite{aufl95}:
\begin{flalign}
\frac{\partial}{\partial x} \arg(V(f,g)(x,\omega)) = \Im \left(\frac{V(f,{g'})(x,\omega)\overline{V(f,g)(x,\omega)}}
{\left|V(f,g)(x,\omega)\right|^2}\right),
\label{eq:DerivFormula}
\end{flalign}
with $g'(t)=\frac{dg}{dt}(t)$. The benefit of this method is that is
does not require unwrapping, instead the phase derivative is computed
by pointwise operations using a second STFT based on the derivative of
the window.
\\

To understand the phase of the STFT more thoroughly, in particular for applications dealing with multipliers, see for example \cite{nathxxl13,balsto09new,uncconv2011}, we conducted related extensive numerical experiments. In the process we observed a rather peculiar phenomenon in the phase derivative, a recurring pattern that appears in a similar way in the neighborhood around any of the zeros of the STFT. The behaviour of the phase derivative close to the singularity always shows the same characteristic shape, i.e., a negative peak followed by a positive one. We describe this phenomenon and provide a complete analytical explanation.   

This paper is organized as follows: In Section \ref{sec:Observations} we report the numerical results. In Section \ref{sec:Analytical} we give a short, instructive, analytical example for this behaviour. In Section \ref{sec:Mathematical} we give the full analytical results.

Results in this paper have partly been reported at a conference \cite{Jaillet2009}, and a preprint of this paper has already been cited in \cite{journals/spl/AugerCF12}.

\section{Numerical Observations}

\label{sec:Observations}

\begin{figure}[hbt]
	\begin{center}
		\includegraphics[width=\figWidth]{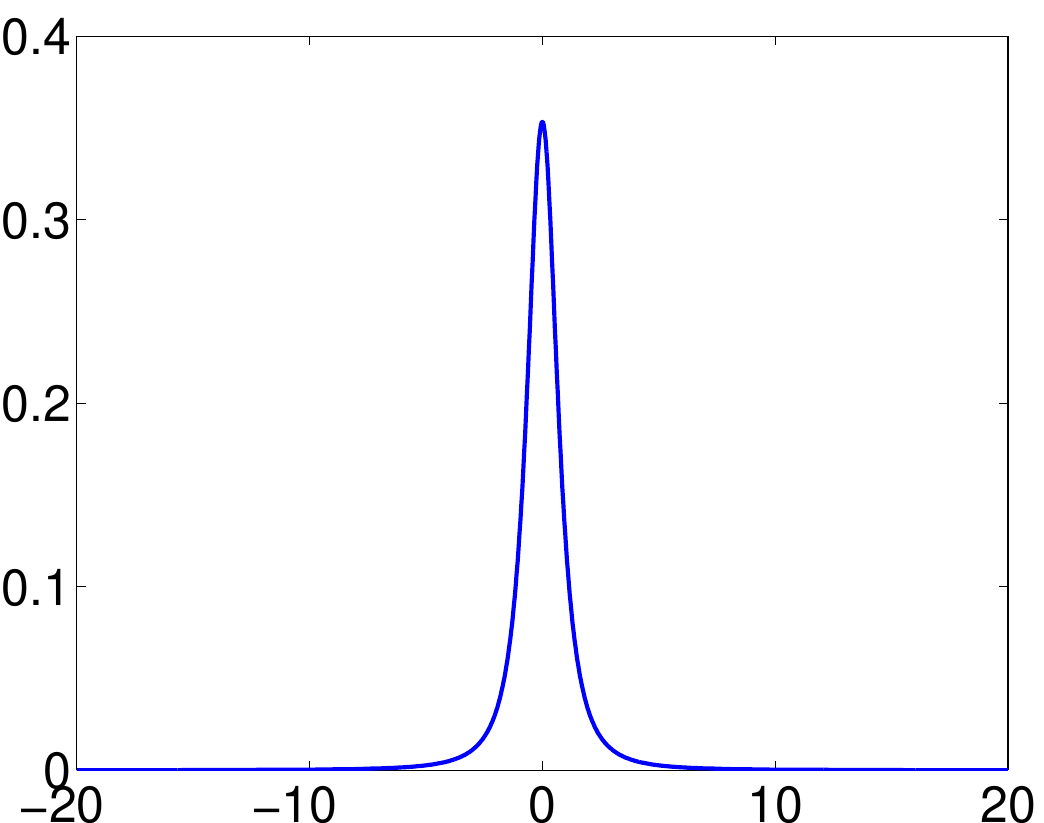}
	\end{center}
	\caption{Distribution of the values of the phase derivative over time of the STFT for a white Gaussian noise.\label{fig:NoiseDistrib}}
\end{figure}

For noise, naturally only statistical properties of the phase are
accessible. Some interesting results for the phase derivative have
been shown in the context of reassignment. In 
\cite{Chassande98on}, the following result is given: We
consider a zero-mean Gaussian analytic white noise $f$ such that
\begin{flalign}
\E[\Re(f(t))\cdot \Re(f(s))]=\E[\Im(f(t))\cdot \Im(f(s))]=\frac{\sigma^2}{2}\delta(t-s)
\end{flalign}
and $\E[f(t) f(s)]=0$ for any $(t,s)\in \bR^2$, with its real and
imaginary parts a Hilbert transform pair. Using a Gaussian window
given by $g(t) = e^{-\pi \frac{t^2}{2\sigma^2}}$, the phase derivative
over time of $V(f,g)$ is a random variable with distribution of the
form: 
\begin{flalign}
\rho(v)=\frac{1}{2(1+v^2)^\frac{3}{2}} .
\end{flalign}
This distribution is shown in Figure \ref{fig:NoiseDistrib}. As can be
seen, it is a quite ``peaky'' distribution, indicating that the values
of the phase derivative are mainly values close to zero, with some
rare values with higher absolute values.

The spatial distribution of the phase derivative seems to be difficult to be solved analytically. Therefore we conducted
systematic numerical experiments to study this spatial distribution.

For this, we need to compute the derivative of the phase in discrete
settings. We used the expression (\ref{eq:DerivFormula}) to compute
the phase derivative.

We see on this formula that we will face numerical difficulties when
the denominator $V(f,g)(x,\omega)$ is close to zero. But using double
precision, these problems only appear for really small values of the
modulus (on the order of $10^{-13}$), which allows us to reliably
observe the values of the phase derivative even close to the zeros of
the STFT. In the figures of this paper, the phase derivative values
are ignored and represented as white at the points where the value of
the modulus is too small. 

\begin{figure}[hbt]
	\begin{center}
		\includegraphics[width=\figWidth]{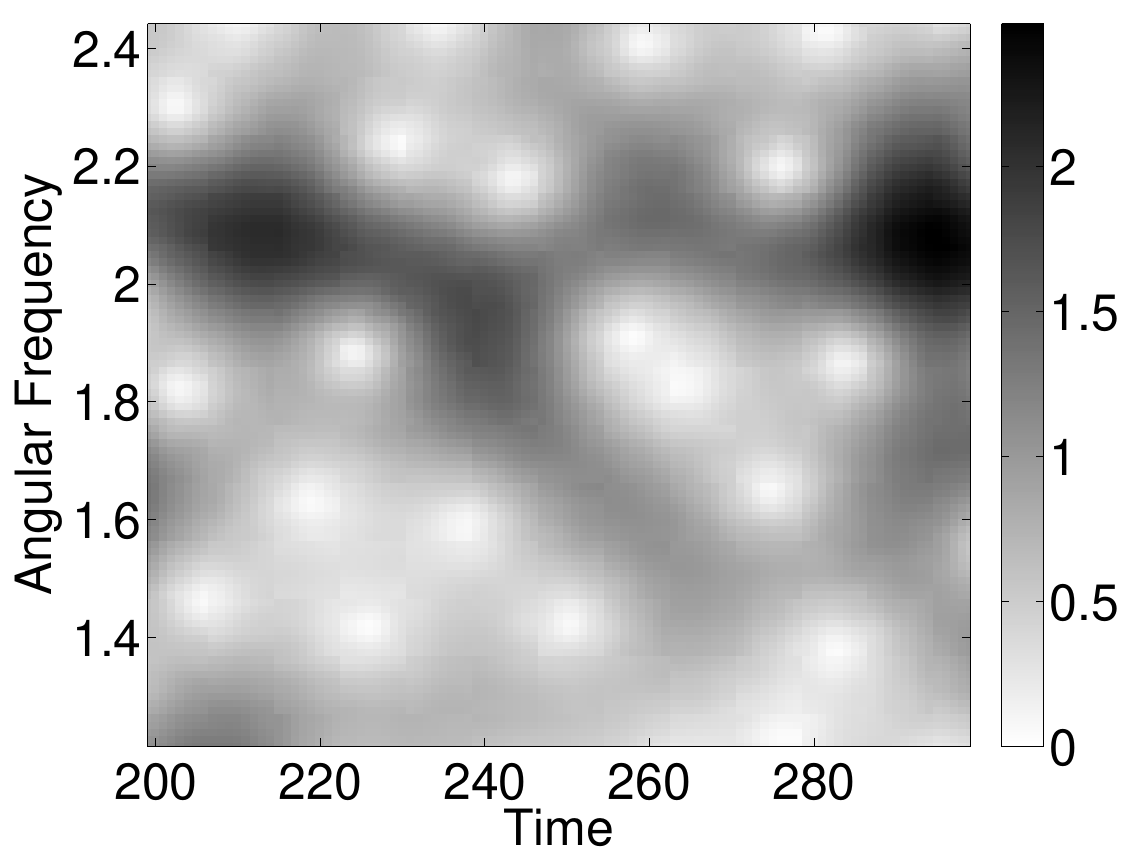}
		
		\includegraphics[width=\figWidth]{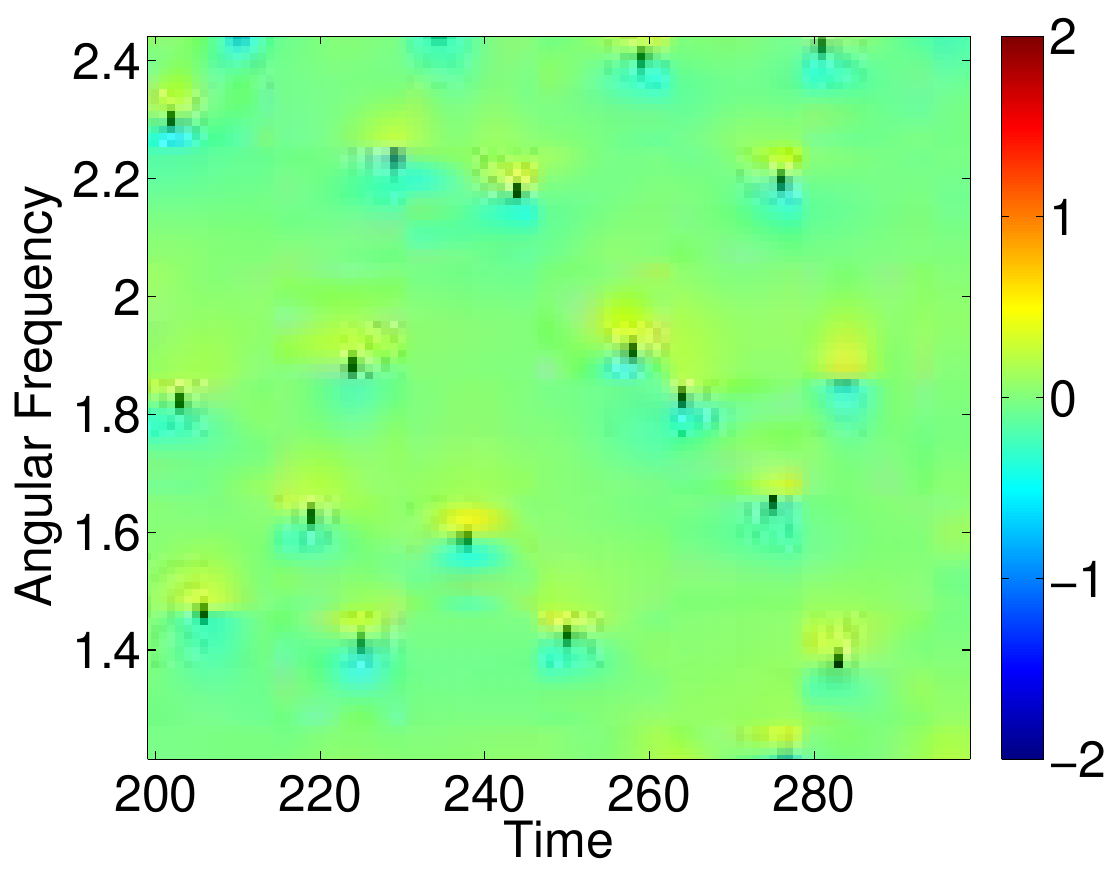}
		\includegraphics[width=\figWidth]{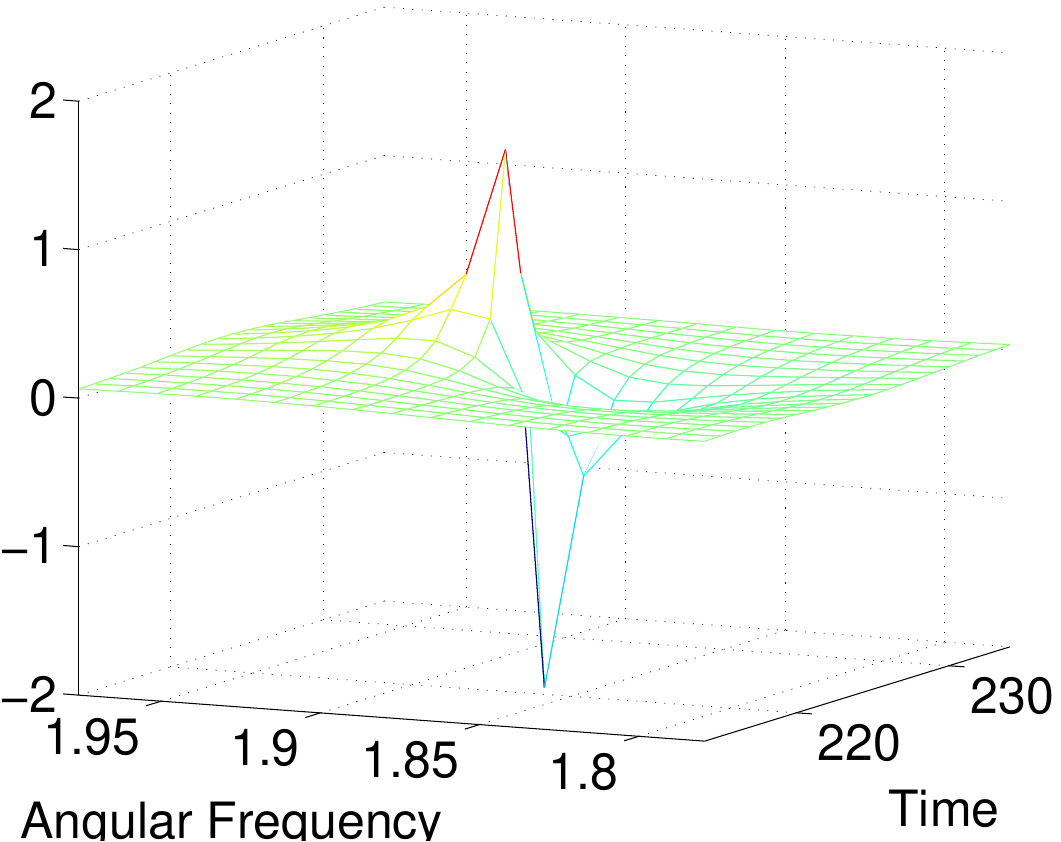}
	\end{center}
	\caption{Observation for a Gaussian white noise, using a Gaussian window. Top: modulus of the STFT. Bottom-left: derivative over time of the phase of the STFT using the definition (\ref{eq:defV}). Bottom-right: mesh plot of the derivative over time of the phase in the neighbourhood of a zero of the STFT.\label{fig:NoiseDerivative}}
\end{figure}

The results of our experiments are illustrated by
Figure \ref{fig:NoiseDerivative}. 
The
time-frequency distribution of the values appears to be highly
structure, as e.g. noted in \cite{gama06}. The values of the phase derivative with
high absolute values are concentrated around several time-frequency
points, which can be identified as the zeros of the transform when
looking at the modulus. Furthermore, the shape of the phase derivative
seems to be very similar in the neighbourhood of the zeros, with a
typical pattern repeating at each zero, see 
Figure
\ref{fig:NoiseDerivative}. When going from low to high frequencies, it
presents a negative peak followed by a positive one.

\begin{figure}[hbt]
	\begin{center}
		\includegraphics[width=\figWidth]{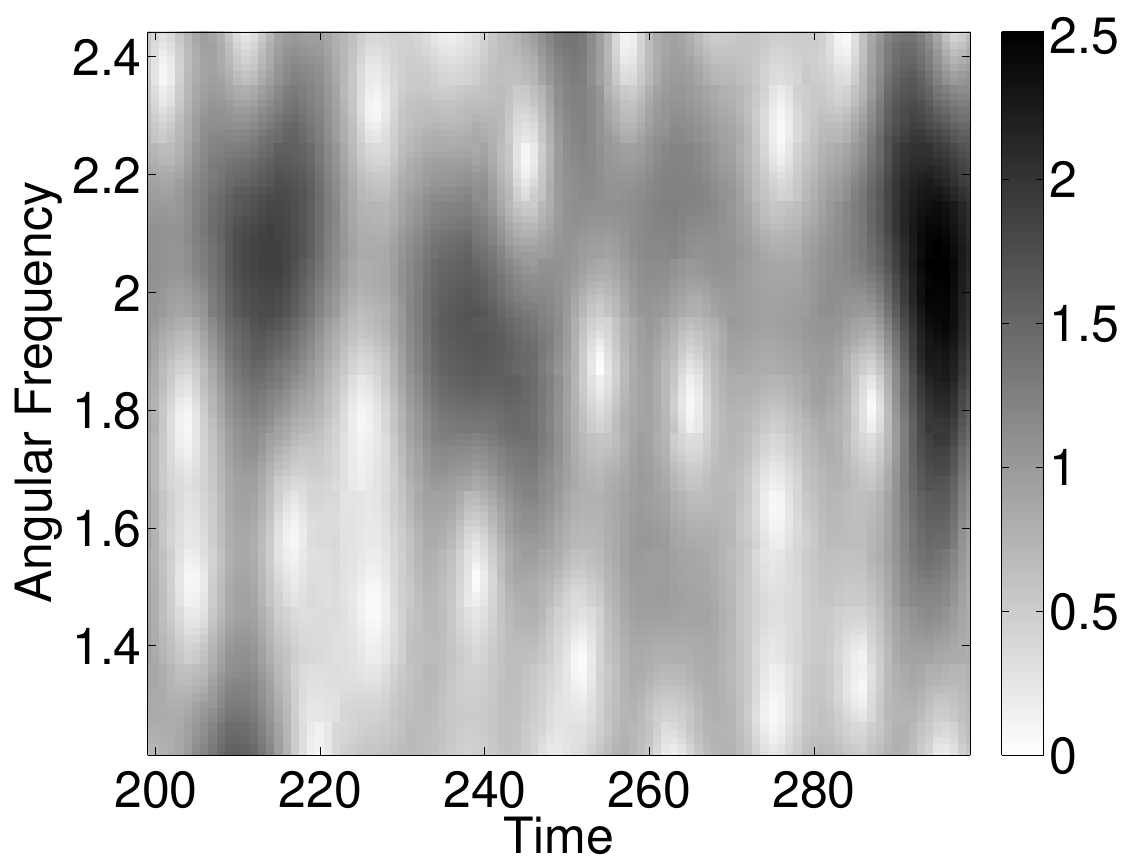}
		\includegraphics[width=\figWidth]{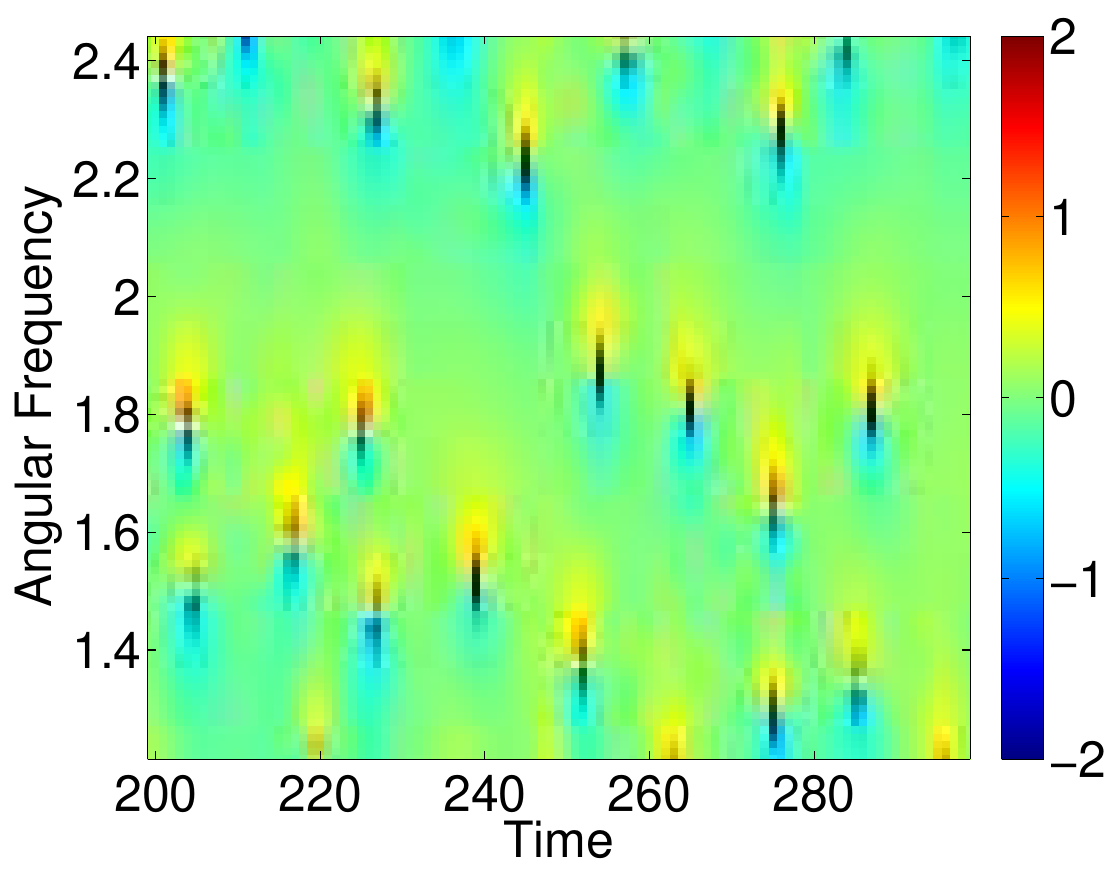}
		\includegraphics[width=\figWidth]{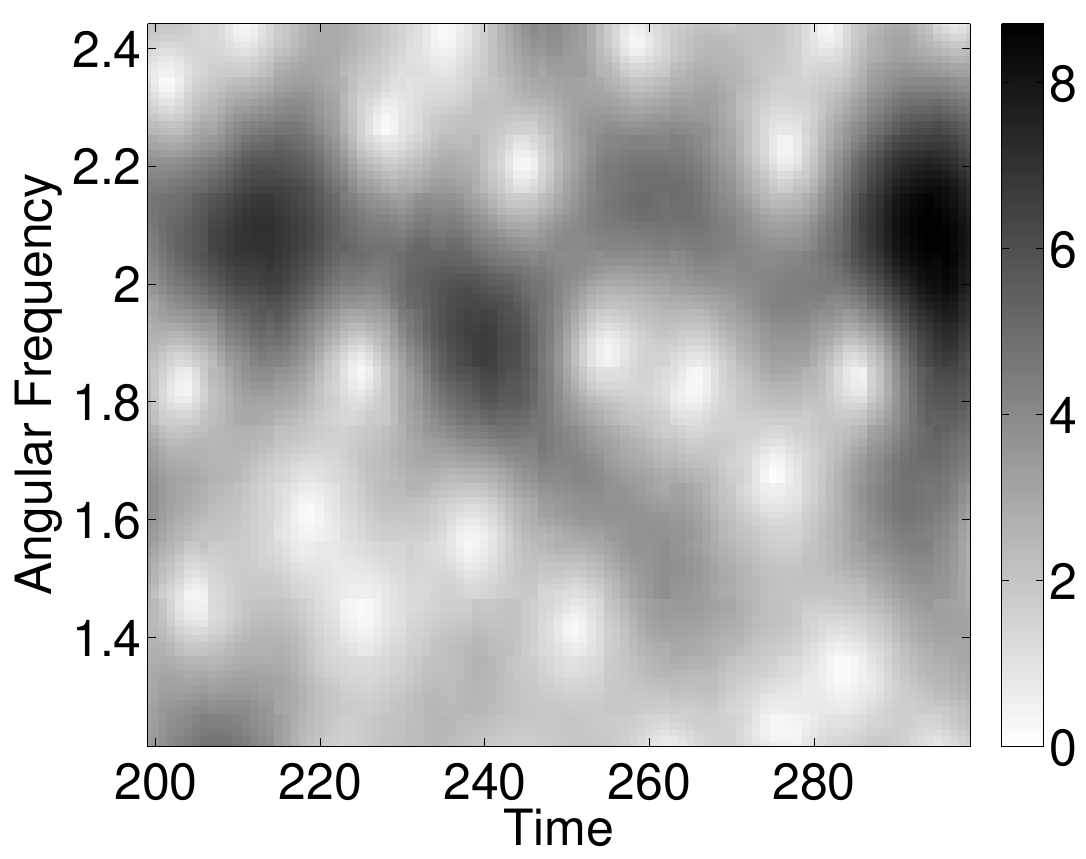}
		\includegraphics[width=\figWidth]{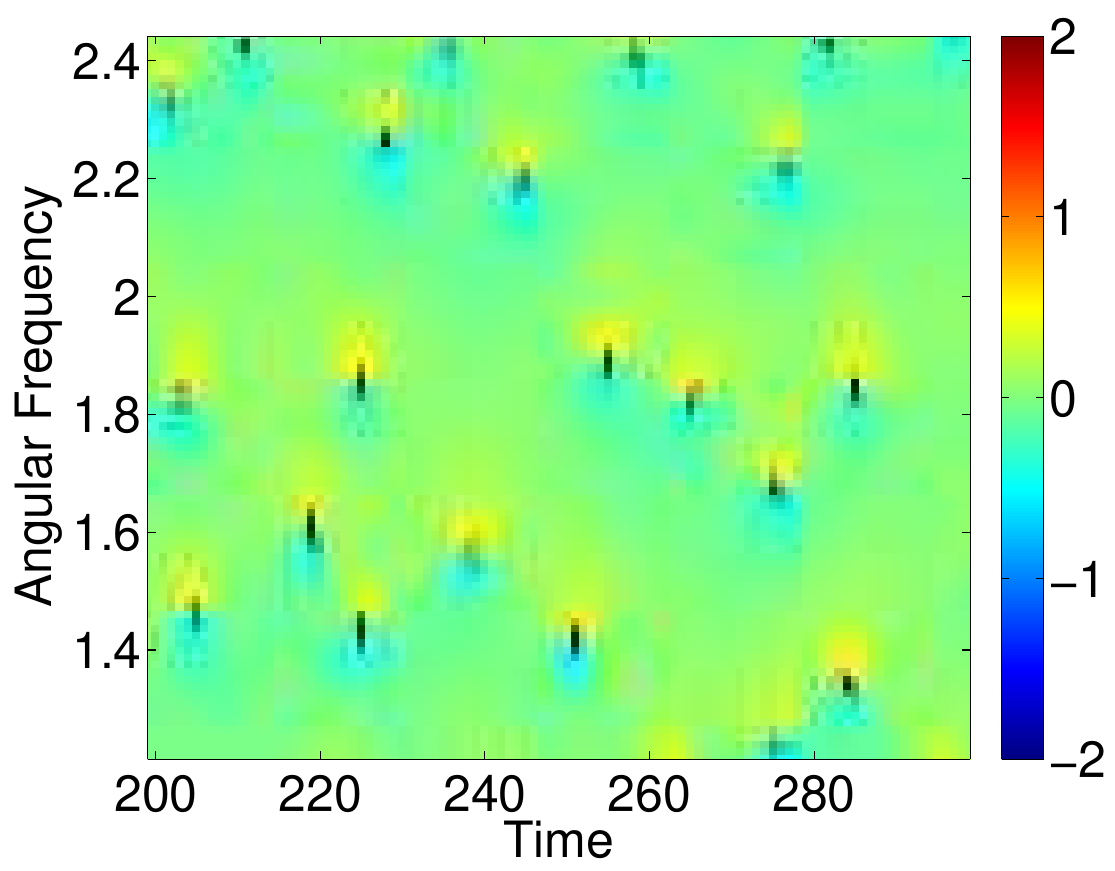}
		\includegraphics[width=\figWidth]{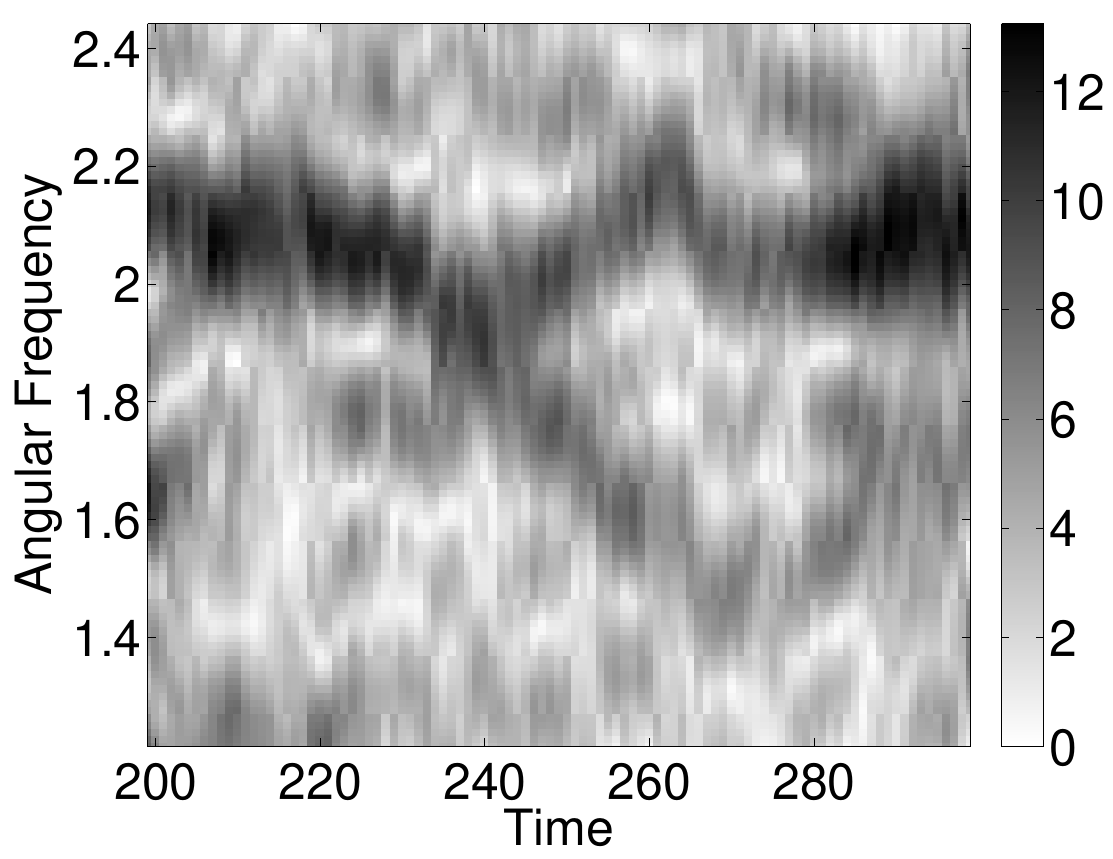}
		\includegraphics[width=\figWidth]{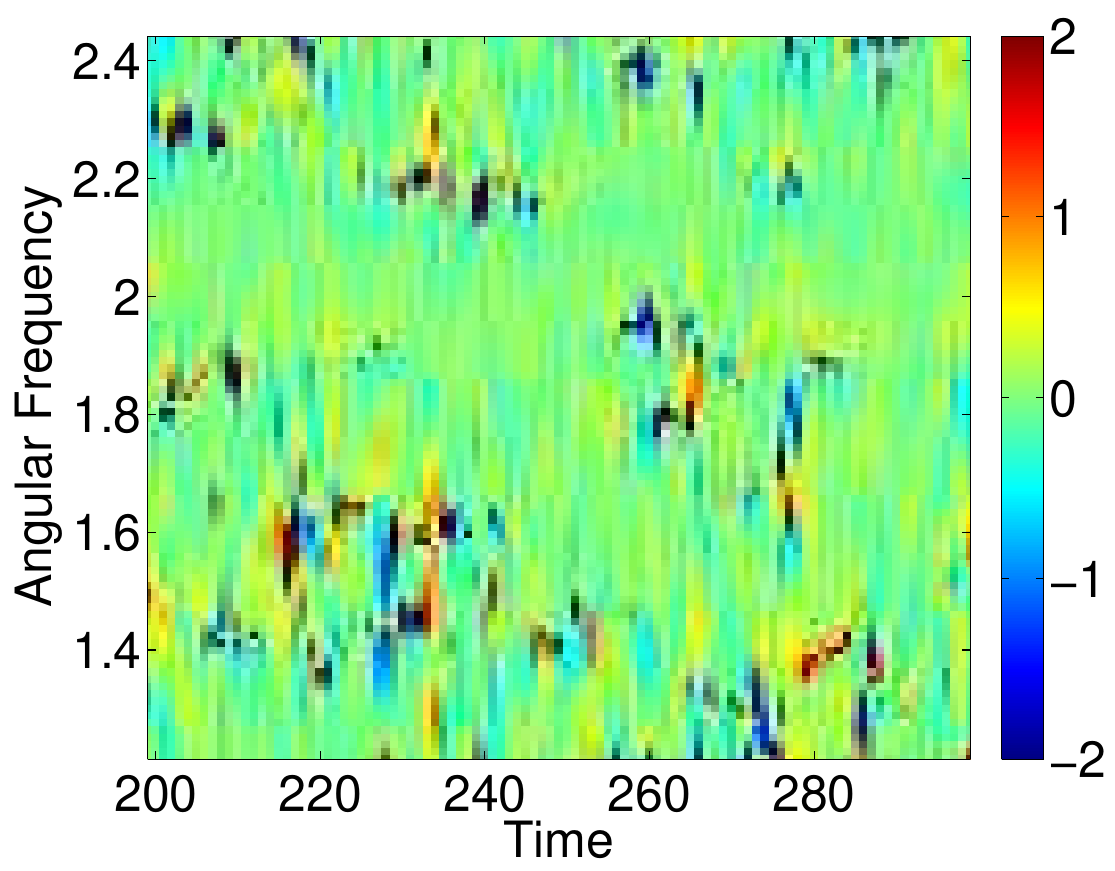}
	\end{center}
	\caption{Influence of the window when analyzing a frozen Gaussian
          white noise. For three different windows, on the left,
          modulus of the STFT, on the right, derivative over time of
          the phase of the STFT using the definition
          (\ref{eq:defV}). From top to bottom, the windows are: a
          narrower Gaussian window, a Hamming window, a rectangular
          window.\label{fig:WinInfluence}}
\end{figure}

This phenomenon is related to the fact that the STFT of white noise is
a correlated process, with a correlation determined by the window
through the reproducing kernel of the transform (see part 6.2.1 of
\cite{carmtorr1}). It is thus interesting to study the influence of
the window choice on the observed structure of the phase derivative, as 
illustrated in Figure
\ref{fig:WinInfluence}. 

We can observe 
that narrowing the window results in similar patterns around the
zeros, but with a scaled shape: the resulting pattern is narrower over
time, but wider over frequency.
Figure \ref{fig:WinInfluence} also shows the influence of the window
type. 
The structure is more
complicated for windows with bad time-frequency concentration. 
On the representation using a Hamming window, we still
observe repeating patterns at the zeros of the transform, but the
variability of the shape of this pattern seems higher, and the pattern
orientation slightly varies, whereas it is fixed in the case of a
Gaussian window.
For the case of the rectangular window, the zeros of the STFT form a more
complicated, extended structures. This leads to much more variable
patterns. Yet, we still, interestingly, observe that the values of the
phase derivative with high absolute values concentrate around the
zeros of the transform, whereas the phase derivative is close to zero
in the regions of the STFT where the modulus is high.

>From the experiments above, we expect those properties to be highly correlated with co-orbit properties regarding the STFT, i.e. inclusion in certain modulation spaces \cite{fe06}. 
In particular, choosing windows in the Feichtinger algebra $S_0$ \cite{feizim1} should result in a phase derivative behavior comparable to the Gaussian window. 
We also expect results to be valid as in Section \ref{sec:Mathematical}, for windows in $S_0$.
The systematic investigation of the phase derivative behaviour for modulation spaces is beyond the scope of this paper, and will be investigate in future work. \\

The behaviour that we observe is not specific to noise signals. Indeed, further experiments on other synthesized and
recorded complex sounds showed that the same characteristics can be
observed for all signals: the values of the phase derivative of high
absolute value are concentrated in the neighbourhood of the zeros of
the STFT, and for ``nice'' windows, a specific pattern appears in this
neighbourhood.

\section{A Simple Explicit Analytic Example}

\label{sec:Analytical}

In this section we give a simple analytical example for which we can explicitly
compute the  phase derivative.

Considering the signal given by
\begin{flalign}
f(t) = e^{2\pi i\omega_1 t}+e^{2\pi i\omega_2 t}
\label{eq:defSig}
\end{flalign}
and using a Gaussian window $g(t) = e^{-\pi \frac{t^2}{2\sigma^2}}$,
we can explicitly compute the expression of the STFT, which results in the formula:
\begin{flalign}
V(f,g)(x,\omega)= e^{-2\pi i x (\omega - \omega_1)}e^{-2\pi \sigma^2 (\omega - \omega_1)^2} + e^{-2\pi i x (\omega - \omega_2)}e^{-2\pi \sigma^2 (\omega -\omega_2)^2}.
\end{flalign}

The zeros of this STFT are the points of coordinates $(x_k, \omega_{mid})$ in the time-frequency plane, 
with $\omega_{mid} = \frac{\omega_1+\omega_2}{2}$ and $ x_k = \frac{1+2k}{2(\omega_1-\omega_2)}$ for $k\in\bZ$.

The expression of the phase derivative for this signal, given in part
VI-12 of \cite{delpal1}, is:
\begin{flalign}
\frac{\partial}{\partial x} \arg(V(f,g)(x,\omega)) = 2\pi \Big( \omega_{mid}- \omega + \delta\tanh(s)\frac{1+\tan^2(2\pi\delta x)}{1+\tan^2(2\pi\delta x)\tanh^2(s)}\Big)
\label{eq:PhaseDeriv}
\end{flalign}
with $\delta=\frac{\omega_2-\omega_1}{2}$ and $s=4\pi \sigma^2 (\omega-\omega_{mid})\delta$.

\begin{figure}[hbt]
	\begin{center}
		\includegraphics[width=\figWidth]{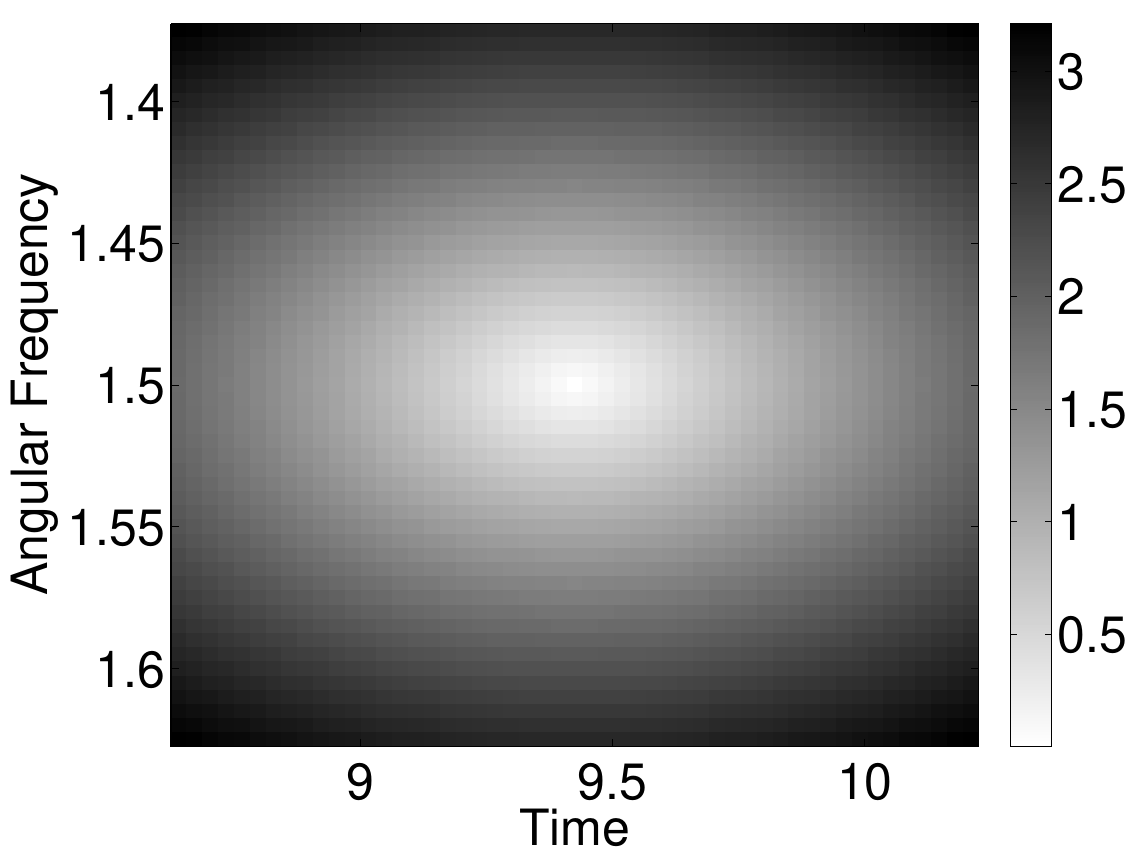}
		
		\includegraphics[width=\figWidth]{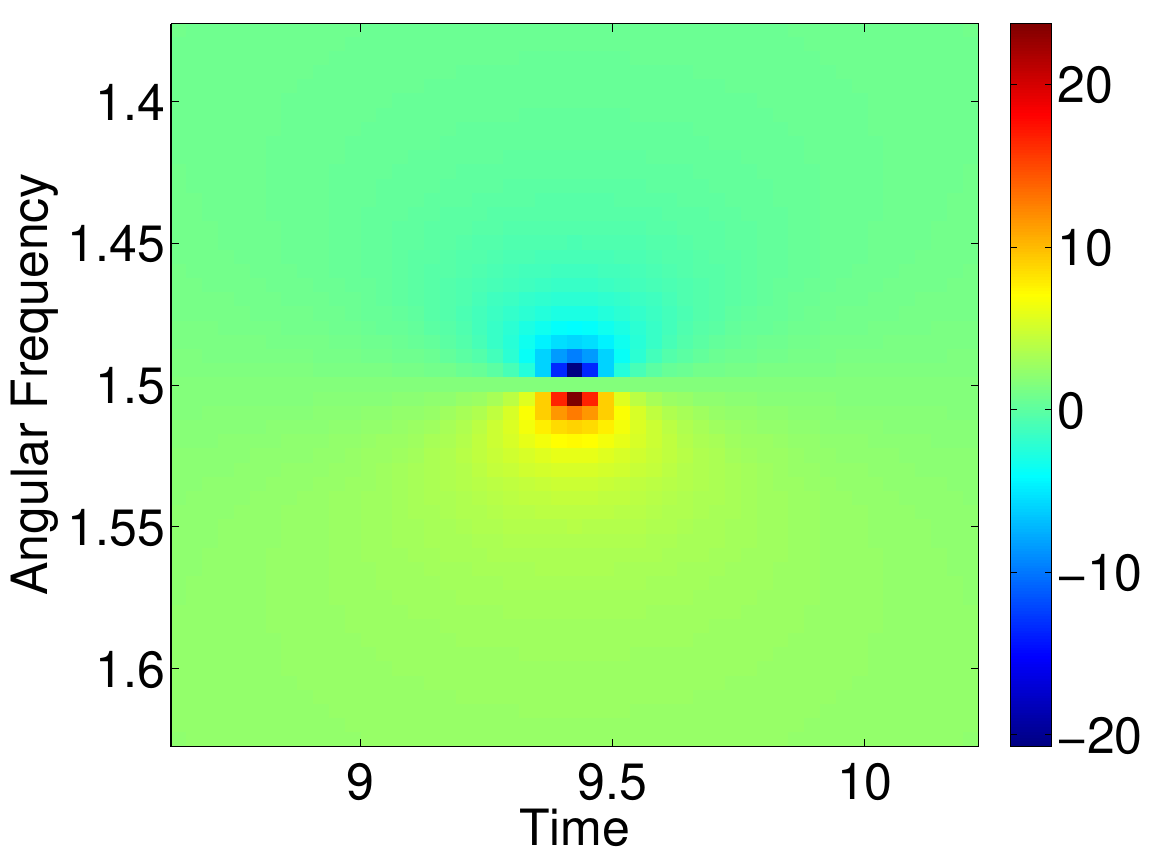}
		\includegraphics[width=\figWidth]{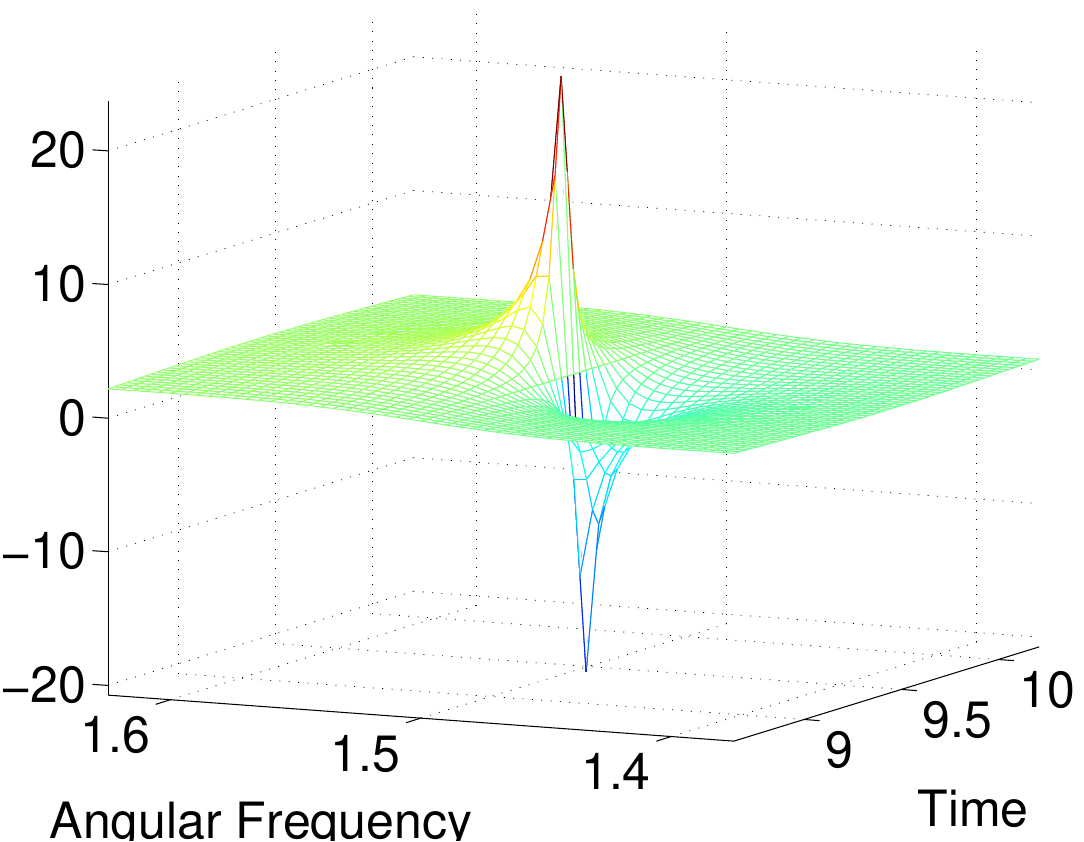}
	\end{center}
	\caption{Observation for the signal defined in
          (\ref{eq:defSig}). Top: modulus of the STFT. Bottom:
          derivative over time of the phase of the STFT according to
          (\ref{eq:PhaseDeriv}) represented as an image (left) and as
          a mesh (right).\label{fig:AnalyticalExamp}}
\end{figure}

The plot of this function around one of the zeros is visible in Figure
\ref{fig:AnalyticalExamp}. We see the pattern that was already observed in the previous section.

For an extended treatment of this analytic example, see \cite{journals/spl/AugerCF12}, where the authors are already referring to a preprint of our present work.

\section{Analytical Results}

\label{sec:Mathematical}

In this section, we denote by $M_h$ the modulation operator $M_h: L^2(\R) \to L^2(\R)$, $f(t) \mapsto M_hf(t) := e^{-2\pi i h t}f(t)$, and by $T_h$ the translation operator $T_h:L^2(\R) \to L^2(\R)$, $f(t) \mapsto T_hf(t) := f(t-h)$ (with $h \in \R$). The set $\swrz(\R)$ is the Schwartz class of rapidly decaying functions. The Fourier transform is denoted by $\mathcal{F}$.

\subsection{Regularity Properties of the STFT}

\label{subsec:Regularity}

\begin{definition}
Define the (unbounded) operator $P$ on $L^2(\R)$ as the multiplication operator
\[
Pf (t) := 2\pi i t\cdot f(t)
\]
with domain 
\[ \mbox{Dom}(P) := \{ f \in L^2(\R): \int_{\R} | t\,f(t) |^2\,dt < \infty \} \subset L^2(\R).
\]
Further, define the
(unbounded) operator
\[
Q := \mathcal{F}^{-1} P \mathcal{F}
\]
(where $\mathcal{F}$ denotes the Fourier transform) with domain 
\[ \mbox{Dom}(Q) := \{ f \in L^2(\R):\mathcal{F}f \in \mbox{Dom}(P) \} \subset L^2(\R).
\]
\end{definition}

In quantum mechanics, these operators are essentially the momentum and position operator, respectively. The operator $P$ (and thus also $Q$) are clearly densely-defined, since $\swrz(\R) \subset \mbox{Dom}(P)$ (and $\mathcal{F}^{-1}\swrz(\R) = \swrz(\R) \subset \mbox{Dom}(Q)$). It can be shown that
$P$ and $Q$ are closed unbounded operators and that $iP$ and $iQ$ are self-adjoint.

We collect basic properties of these operators in the following lemma.

\begin{lemma}
The operators $P$ and $Q$ have the following properties:
\begin{enumerate}
\item $\mathcal{F}Q = P\mathcal{F}$ on $\mbox{Dom}(Q)$, $Q\mathcal{F} = - \mathcal{F}P$ on $\mbox{Dom}(P)$;
\item $Q$ is a (maximal extension of a) differential operator, more precisely: if $f \in \swrz(\R)$, then $Qf(t) = \frac{d}{dt} f(t)$.
\end{enumerate}
\end{lemma}

The next lemma is in essence a well-known result from the theory of operator (semi-)groups; it gives the infinitesimal generators of the modulation and translation group, respectively, see e.g. \cite{fosi97,gr01}. 
The version below needed in this manuscript can be proved in a straight-forward way.

\begin{lemma}
Let $f \in \mbox{Dom}(P) $, then
\[
\| \frac{1}{h}(M_h - Id)f - Pf \|_{L^2} \to 0 
\]
for $h \to 0$.\\
Let $f \in \mbox{Dom}(Q)$, then
\[
\| \frac{1}{h}(T_h - Id)f + Qf \|_{L^2} \to 0 
\]
for $h \to 0$.
\end{lemma}

We can now prove a regularity result for the short-time Fourier transform.

\begin{prop}\label{T:main_regularity_thm}
Let $f,g \in L^2(\R)$.
\begin{enumerate}
\item If $f$ belongs to $\mbox{Dom}(P)$, then $V(f,g)$ has a continuous partial derivative with respect to the second argument $\omega$, and we have
\[
\frac{\partial}{\partial \omega} V(f,g)(x, \omega) = - V(Pf,g)(x, \omega).
\]
\item If $g$ belongs to $\mbox{Dom}(Q)$, then $V(f,g)$ has a continuous partial derivative with respect to the first argument $x$, and we have
\[
\frac{\partial}{\partial x} V(f,g)(x, \omega) = - V(f, Qg)(x, \omega).
\]
\end{enumerate}
\end{prop}
\begin{proof}
Assume $f \in \mbox{Dom}(P)$. Then, by the preceding lemma,
\begin{align*}
\frac{1}{h}(V(f,g)(x, \omega + h) & - V(f,g)(x,\omega)) = \langle f, \frac{M_h - Id}{h} M_{\omega}T_x g \rangle \\
& = \langle \frac{M_{-h} - Id}{h}f,  M_{\omega}T_x g \rangle \\
& \stackrel{h \to 0}{\longrightarrow} \langle -Pf,  M_{\omega}T_x g \rangle = - V(Pf,g)(x,\omega),
\end{align*}
which is a continuous function on $\R^2$.\\
If $g \in \mbox{Dom}(Q)$, then, analogously,
\begin{align*}
\frac{1}{h}(V(f,g)(x + h, \omega) & - V(f,g)(x,\omega)) = \langle f, M_{\omega}T_x \frac{T_h - Id}{h} g \rangle \\
& \stackrel{h \to 0}{\longrightarrow} \langle f,  - M_{\omega}T_x Qg \rangle = - V(f,Qg)(x,\omega).
\end{align*}
\end{proof}

Using $V(f,g)(x, \omega) = e^{-2\pi i x \omega} \overline{V(g,f)(-x, -\omega)}$, the partial derivatives of $V(f,g)$ exist if and only if those of $V(g,f)$ exist. We may thus change the roles of $f$ and $g$.
\begin{corollary}\label{C:cor_main_regularity_thm}
Let $f,g \in L^2(\R)$.
\begin{enumerate}
\item If $g$ belongs to $\mbox{Dom}(P)$, then $V(f,g)$ has a continuous partial derivative with respect to the second argument $\omega$, and we have
\[
\frac{\partial}{\partial \omega} V(f,g)(x, \omega) = V(f,Pg)(x, \omega) - 2\pi i x \,V(f,g)(x,\omega).
\]
\item If $f$ belongs to $\mbox{Dom}(Q)$, then $V(f,g)$ has a continuous partial derivative with respect to the first argument $x$, and we have
\[
\frac{\partial}{\partial x} V(f,g)(x, \omega) = V(Qf, g)(x, \omega) - 2 \pi i \omega\, V(f,g)(x,\omega).
\]
\end{enumerate}
\end{corollary}

If $f$ (resp. $g$) belongs to Schwartz class $\swrz(\R)$, then $Pf, Qf$ (resp. $Pg, Qg$) $\in \swrz(\R)$, as well. Iterated application of Proposition \ref{T:main_regularity_thm} or Corollary \ref{C:cor_main_regularity_thm} gives the following smoothness result for the STFT.

\begin{theorem}
Let $f,g \in L^2(\R)$, and at least one of them in Schwartz class $\swrz(\R)$. Then $V(f,g)(x,\omega)$ is infinitely partially differentiable in both variables $x$ and $\omega$.
\end{theorem}

Although this result may be considered mathematical folklore, to our knowledge it has not been stated and proved in the literature so far.
Note that this proves in particular that the STFT with Gaussian window is smooth.

\subsection{The Derivative of the Phase Around the Zeros of the STFT}

\label{subsec:Derivative}

In this section we present an analytic explanation of the peculiar
behaviour of the phase derivatives of the STFT for a large class of
window functions. It turns out that the phenomenon is connected to the
smoothness and continuous differentiability of the STFT which, as we
have seen in the previous paragraph, is in turn connected to the smoothness of the window.

Consider first the partial derivative of the phase of the STFT with
respect to the first variable (i.e., the 'time' variable). For convergence along a
vertical path, we have

\begin{theorem}[\label{thm:derivI}Phase derivatives of the STFT, part I]
Let $f, g \in L^2(\R)$. Assume that
\begin{itemize}
\item $V(f,g) = V = U + i\cdot W \in C^2(\R^2)$
\item $V(x_0, \omega_0) = 0$
\item $\det{J_V}(x_0, \omega_0) < 0$, where
\[
J_V(x_0, \omega_0) = \left(\begin{matrix}
\displaystyle \frac{\partial U}{\partial x}{(x_0, \omega_0)} & \displaystyle \frac{\partial U}{\partial \omega}{(x_0, \omega_0)} \\ &  \\ \displaystyle \frac{\partial W}{\partial x}{(x_0, \omega_0)} & \displaystyle \frac{\partial W}{\partial \omega}{(x_0, \omega_0)}
\end{matrix}\right)
\]
denotes the Jacobian matrix of $V$ at the point $(x_0, \omega_0)$
\end{itemize}
Then the phase $\psi(x,\omega)$ of $V(f,g)(x,\omega)$ satisfies
\[
\frac{\partial \psi}{\partial x}(x_0,\omega) \longrightarrow 
\begin{cases}
+\infty, &\text{if $\omega \uparrow \omega_0$ from below}\\
-\infty, &\text{if $\omega \downarrow \omega_0$ from above}.
\end{cases}
\]
\end{theorem}
\begin{proof}
We have
\[
\frac{\partial \psi}{\partial x}(x_0, \omega) =  \frac{U(x_0,\omega)\cdot W_x(x_0,\omega) - W(x_0,\omega) \cdot U_x(x_0,\omega)}{U^2(x_0,\omega) + W^2(x_0,\omega)}.
\]
Since $W_x$ and $U_x$ are continuous and thus remain bounded in a neighbourhood of $(x_0, \omega_0)$, both numerator and denominator tend to zero for $\omega \to \omega_0$. However, both functions are differentiable, since $V \in C^2(\R^2, \R^2)$. So L'Hospital's Rule is applicable and yields
\begin{align*}
\lim_{\omega\to\omega_0} & \frac{U(x_0,\omega)\cdot W_x(x_0,\omega) - W(x_0,\omega) \cdot U_x(x_0,\omega)}{U^2(x_0,\omega) + W^2(x_0,\omega)} \\
& \stackrel{\mbox{\tiny L'Hosp.}}{=} \lim_{\omega\to\omega_0} \frac{\frac{d}{d\omega}\left(U(x_0,\omega)\cdot W_x(x_0,\omega) - W(x_0,\omega) \cdot U_x(x_0,\omega)\right)}{\frac{d}{d\omega}\left(U^2(x_0,\omega) + W^2(x_0,\omega)\right)} \\
& = \lim_{\omega\to\omega_0} \frac{\left(U_{\omega} W_x + U W_{x\omega} - W_{\omega} U_x - W U_{x\omega}\right)(x_0,\omega)}{\left(2 U U_{\omega} + 2 W W_{\omega}\right)(x_0,\omega)} \\
& = \lim_{\omega\to\omega_0} \frac{\left(U W_{x\omega} - W U_{x\omega}\right)(x_0,\omega) + \left( U_{\omega} W_{x} - W_{\omega} U_{x} \right)(x_0,\omega)}{\left(2 U U_{\omega} + 2 W W_{\omega}\right)(x_0,\omega)}.
\end{align*}
Concerning this limit, we clearly have
\[
\left(U W_{x\omega} - W U_{x\omega}\right)(x_0,\omega) \to 0,
\]
since $U(x_0,\omega) \to U(x_0,\omega_0) = 0$, $W(x_0,\omega) \to W(x_0,\omega_0) = 0$, and $W_{x\omega}$ and $U_{x\omega}$ are continuous and thus remain bounded in a neighborhood of $(x_0, \omega_0)$. Furthermore,
\begin{align*}
\left( U_{\omega} W_{x} - W_{\omega} U_{x} \right)(x_0,\omega) & \to \left( U_{\omega} W_{x} - W_{\omega} U_{x} \right)(x_0,\omega_0) \\
& \quad \quad = - \det{J_V}(x_0,\omega_0) \neq 0,
\end{align*}
by assumption. Hence the numerator tends to a nonzero number, in this case ($\det{J_V}(x_0, \omega_0) < 0$) a positive one. For the denominator, we find
\begin{align*}
\left(2 U U_{\omega} + 2 W W_{\omega}\right)(x_0,\omega) & = \frac{d}{d\omega}\left(U^2 + W^2\right)(x_0,\omega) \\
& = \begin{cases}
< 0, & \text{ if $\omega < \omega_0$} \\
> 0, & \text{ if $\omega > \omega_0$}
\end{cases}
\end{align*} 
since the function $\omega \mapsto \left(2 U U_{\omega} + 2 W W_{\omega}\right)(x_0,\omega)$ has a strict local minimum in $\omega_0$. At the same time,
\[
\left(2 U U_{\omega} + 2 W W_{\omega}\right)(x_0,\omega) \to 0
\]
for $\omega \to \omega_0$, hence the denominator goes to zero from below for $\omega \uparrow \omega_0$ and from above for $\omega \downarrow \omega_0$. This concludes the proof.
\end{proof}

Note that for simplicity we have only considered the case that the Jacobian determinant $\det{J_V}(x_0, \omega_0)$ is negative; this case corresponds to the examples we presented above. For positive Jacobian determinant, the situation is completely analogous, although reversed in the sense that the positive and negative singularities switch roles. Apart from this, the general behaviour remains the same.

For convergence along a horizontal path, we need  slightly more regularity:

\begin{theorem}[\label{thm:derivII}Phase derivatives of the STFT, part II]
Let $f, g \in L^2(\R)$. Assume that
\begin{itemize}
\item $V(f,g) = V = U + i\cdot W \in C^3(\R^2)$
\item $V(x_0, \omega_0) = 0$
\item $\det{J_V}(x_0, \omega_0) < 0$, where $J$ is the Jacobian as in Theorem \ref{thm:derivI}
\end{itemize}
Then there exists a number $c \in \R$ such that the phase $\psi(x,\omega)$ of $V(f,g)(x,\omega)$ satisfies
\[ \lim_{x\to x_0} \frac{\partial \psi}{\partial x}(x,\omega_0) = c.
\]
\end{theorem}
\begin{proof}
The assumptions allow us to apply L'Hospital's Rule twice, giving
\begin{align*}
\lim_{x\to x_0} & \frac{\partial \psi}{\partial x}(x,\omega_0) \\
& = \lim_{x\to x_0} \frac{U(x,\omega_0)\cdot W_x(x,\omega_0) - W(x,\omega_0) \cdot U_x(x,\omega_0)}{U^2(x,\omega_0) + W^2(x,\omega_0)} \\
& \stackrel{\mbox{\tiny L'Hosp.}}{=} \lim_{x\to x_0} \frac{\frac{d}{dx}\left(U(x,\omega_0)\cdot W_x(x,\omega_0) - W(x,\omega_0) \cdot U_x(x,\omega_0)\right)}{\frac{d}{dx}\left(U^2(x,\omega_0) + W^2(x,\omega_0)\right)} \\
& = \lim_{x\to x_0} \frac{\left(U_x W_x + U W_{xx} - W_x U_x - W U_{xx}\right)(x,\omega_0)}{\left(2 U U_{x} + 2 W W_{x}\right)(x,\omega_0)} \\
& = \lim_{x\to x_0} \frac{\left(U W_{xx} - W U_{xx}\right)(x,\omega_0)}{\left(2 U U_{x} + 2 W W_{x}\right)(x,\omega_0)} \\
& \stackrel{\mbox{\tiny L'Hosp.}}{=} \lim_{x\to x_0} \frac{\frac{d}{dx}\left(U W_{xx} - W U_{xx}\right)(x,\omega_0)}{\frac{d}{dx}\left(2 U U_{x} + 2 W W_{x}\right)(x,\omega_0)} \\
& = \lim_{x\to x_0} \frac{\left( U_xW_{xx} + UW_{xxx} - W_xU_{xx} - WU_{xxx} \right)(x,\omega_0)}{2\left( U_x^2 + UU_{xx} + W_x^2 + WW_{xx} \right)(x,\omega_0)}.
\end{align*}
For the denominator,
\[
\left(UU_{xx} + WW_{xx}\right)(x,\omega_0) \to \left(UU_{xx} + WW_{xx}\right)(x_0,\omega_0) = 0
\]
for $x\to x_0$, but
\[
\left( U_x^2 + W_x^2 \right)(x,\omega_0) \to \left( U_x^2 + W_x^2 \right)(x_0,\omega_0) > 0
\]
converges to a nonzero number, since not both $U_x(x_0,\omega_0)$ and
$W_x(x_0, \omega_0$ can be zero because of $\det{J_V}(x_0,\omega_0) =
\left( U_xW_{\omega} - U_{\omega}W_x\right)(x_0, \omega_0) \neq
0$. The numerator obviously converges:
\begin{align*}
& \left( U_xW_{xx} + UW_{xxx} - W_xU_{xx} - WU_{xxx} \right)(x,\omega_0) \\
& \quad \quad \to \left( U_xW_{xx} - W_xU_{xx} \right)(x_0,\omega_0) \in \R,
\end{align*}
thus
\[
\lim_{x\to x_0} \frac{\partial \psi}{\partial x}(x,\omega_0) = \frac{\left( U_xW_{xx} - W_xU_{xx} \right)(x_0,\omega_0)}{2\left( U_x^2 + W_x^2 \right)(x_0,\omega_0)} =: c \in \R.
\]
\end{proof}

Concerning the partial derivatives of the phase of the STFT with
respect to the second variable (i.e., the 'frequency' variable), we can argue almost
identically and thus find the following analogous results:

\begin{theorem}[\label{thm:derivIII}Phase derivatives of the STFT, part III]
Let $f, g \in L^2(\R)$. Assume that
\begin{itemize}
\item $V(x_0, \omega_0) = 0$
\item $\det{J_V}(x_0, \omega_0) < 0$
\end{itemize}
Let $V(f,g) = V = U + i\cdot W \in C^2(\R^2)$.
Then the phase $\psi(x,\omega)$ of $V(f,g)(x,\omega)$ satisfies
\[
\lim_{x \to x_0} \frac{\partial \psi}{\partial \omega}(x,\omega_0) = 
\begin{cases}
+\infty, &\text{if $x \rightarrow x_0$ from the left}\\
-\infty, &\text{if $x_0 \leftarrow x$ from the right}.
\end{cases}
\]
Let $V(f,g) = V = U + i\cdot W \in C^3(\R^2)$, then
\[ \lim_{\omega\to \omega_0} \frac{\partial \psi}{\partial \omega}(x_0,\omega) = c^{\prime} \in \R, \quad \text{if $\omega \to \omega_0$},
\]
converges to some real number $c^{\prime} \in\R$. \hfill$\Box$
\end{theorem}

\section*{Acknowledgements}

\label{sec:Acknowledgement}

The work on this paper was partly supported by the Austrian Science Fund (FWF) START-project FLAME ('Frames and Linear Operators for Acoustical Modeling and Parameter Estimation'; Y 551-N13) and the Vienna Science and Technology Fund (WWTF) project CHARMED ('Computational harmonic analysis of high-dimensional biomedical data'; VRG12-009). 
The authors are
thankful to the project partners for fruitful discussions and
valuable comments, in particular to B. Torr\'esani and M. Ehler. 
Particular thanks go to P. Flandrin for his kind encouragement. 
The authors would also like to thank the anonymous reviewer for valuable comments!

\end{document}